\documentclass{article}
\usepackage{latexsym,amsfonts,amssymb}
\usepackage[latin1]{inputenc}
\usepackage{amssymb}
\usepackage{mathrsfs}
\usepackage{float}
\usepackage{amsmath,amsfonts,amssymb,amsthm}
\usepackage{caption}
\usepackage{enumitem}

\theoremstyle{definition}

\title{\sc An Arad and Fisman's theorem on products of conjugacy classes revisited}

\author{
 Antonio Beltr\'an\\
\footnotesize
Departamento de Matem\'aticas,\\
\footnotesize Universidad Jaume I, \footnotesize
12071 Castell\'on, Spain\\
\footnotesize
e-mail: abeltran@uji.es\\
\\
 Mar\'{\i}a Jos\'e Felipe\\
\footnotesize
Instituto de Matem\'atica Pura y Aplicada,\\
\footnotesize Universitat Polit\`ecnica de Val\`encia,  \footnotesize
46022, Valencia, Spain,\\
\footnotesize
e-mail: mfelipe@mat.upv.es\\
\\
Carmen Melchor\\
\footnotesize Departamento de Did\'actica de la Matem\'atica,\\
\footnotesize Universitat de Val\`encia,
\footnotesize Valencia, Spain\\
\footnotesize
e-mail: Carmen.Melchor-Borja@uv.es}

\date{}
\begin{document}
\maketitle

\thanks{}
\begin{abstract}
 A theorem of Z. Arad and E. Fisman establishes that if $A$ and $B$ are two conjugacy classes of a  finite group $G$ such that either $AB=A\cup B$ or $AB=A^{-1} \cup B$, then $G$ cannot be non-abelian simple. We demonstrate that, in fact, $\langle A\rangle = \langle B\rangle$ is solvable, the elements of $A$ and $B$ are $p$-elements for some prime $p$, and $\langle A\rangle $ is $p$-nilpotent. Moreover, under the second assumption, it turns out that $A=B$ and this is the only possible case. This  research is done by appealing to recently developed techniques and results that are based on the Classification of Finite Simple Groups.
 \end{abstract}



\section{Introduction}

The well-known and long-standing conjecture of Arad and Fisman claims that the product of two non-trivial conjugacy classes of a finite non-abelian simple group cannot be a conjugacy class.  Taking one further step, several authors have studied more general conditions on  the product of conjugacy classes that cannot either happen in a non-abelian simple group. This occurs, for instance, when the product of two conjugacy classes is the union of some restricted sets of conjugacy classes (see for instance \cite{AradFismanMuzychuk}, \cite{AradMuzychuk}, \cite{Nuestro6}).

\bigskip
Our contribution is motivated by the new techniques and results that have been developed in the last few years  (requiring the Classification of Finite Simple Groups) for making progress in this direction, and these not only provide the non-simplicity of the group but the solvability of certain subgroups generated by the conjugacy classes under certain conditions on their products. This is the case, for example, of the main results of \cite{GuralnickNavarro}, \cite{Powers} and \cite{Nuestro5}.

\bigskip
 Among these results, Arad and Fisman  proved that  when $A$ and $B$ are two conjugacy classes of a group $G$ such that either $AB=A\cup B$ or $AB=A^{-1} \cup B$, then $G$ cannot be non-abelian simple \cite{AradFisman}. They used elementary methods to prove this, however, the new outlined approaches allow us to revisit this theorem and supply solvability  and structural properties within the group. We prove the following.

\bigskip
\noindent
\textbf{Theorem A.} {\it Let $A$ and $B$ be conjugacy classes of a finite group $G$ and suppose that $AB= A \cup B$. Then $\langle A\rangle= \langle B\rangle$ is solvable. Furthermore, the elements of $A$ and $B$ are $p$-elements for some prime $p$ and $\langle A\rangle$ is $p$-nilpotent.}

\bigskip
\noindent
\textbf{Theorem B.} {\it  Let $A$ and $B$ be conjugacy classes of a finite group $G$ and suppose that $AB= A^{-1} \cup B$ with $A\neq A^{-1}$. Then $A=B$ and $\langle A\rangle$ is solvable. Furthermore, the elements of $A$ are $p$-elements for some prime $p$ and $\langle A\rangle$ is $p$-nilpotent. }

\bigskip
The case $A=B$ in Theorem B, which means that $A^2=A \cup A^{-1}$, was already studied in  Theorem D of \cite{Powers}. This asserts that, under this hypothesis, $\langle A\rangle$ is solvable and the elements of $A$ are $p$-elements. We will improve this result by showing that $\langle A\rangle$ is in addition $p$-nilpotent. Moreover, for proving Theorem B, we also need a new solvability criterion concerning the product of a conjugacy class and its inverse class, which has interest on its own.

\bigskip 
\textbf{Theorem C.} \textit{Let $A$ be a conjugacy class of a finite group} $G$ \textit{such that}
$AA^{-1}=1 \cup A\cup A^{-1}$. \textit{Then} $\langle A\rangle = AA^{-1}$ \textit{is an elementary abelian group.}

\bigskip
The above result provides further evidence of the following conjecture posed in 
\cite{Nuestro5}. If $A$ and $B$ are conjugacy classes of a group such that $AA^{-1}=1 \cup B\cup B^{-1}$, then $\langle A\rangle$ is solvable. The non-simplicity of $G$ and the solvability of $\langle A\rangle$ in some particular cases were obtained in Theorems A and C of \cite{Nuestro6}. 

\bigskip
The proofs of Theorems A and B are based on the Classification. However, the proof of Theorem C is elementary. All groups are supposed to be finite and the notation is standard and essentially follows the one appearing in \cite{AradFisman}.

\section{Preliminary results}
We state some preliminary results. The first one is necessary for proving both Theorems A and B and part (a) requires the Classification of Finite Simple Groups. However, part (b) does not need it.\\

\noindent
\textbf{Theorem 2.1.} {\it Let $G$ be a finite group and let $N$ be a normal subgroup of $G$. Let
$x\in G$ be such that all elements of $xN$ are conjugate in G. Then:
\begin{enumerate}
\item [{\rm (a)}] $N$ is solvable.
\item [{\rm (b)}] If $x$ is a $p$-element for some prime $p$, then $N$ has a normal $p$-complement.
\end{enumerate}}

{\it Proof.} This is Theorem 3.2 (a) and (c) of \cite{GuralnickNavarro}. $\Box$\\

The following property, however, is elementary and is used for proving Theorem B. Observe that in the particular case of Theorem C, this property is trivial. This situation is addressed in \cite{Nuestro6}.\\

\noindent
\textbf{Lemma 2.2.} {\it Let $K$ and $D$ conjugacy classes of a finite group $G$ such that $KK^{-1}=1\cup D\cup D^{-1}$. If $K$ is real, then $D$ is real.}\\

{\it Proof.} See Lemma 3.1 of \cite{Nuestro6}. $\Box$\\

Our approach mainly utilizes the complex group algebra. We denote by $\mathbb{C}[G]$ the complex group algebra of a group $G$ over the complex field $\mathbb{C}$. Let $K$ be a conjugacy class of $G$ and denote by $\widehat{K}$ the class sum of the elements of $K$ in $\mathbb{C}[G]$. Let $g_1, \ldots, g_k$ be representatives of the conjugacy classes of a finite group $G$. Let $\widehat{S}=\sum_{i=1}^{k}n_i\widehat{g_i^G}$ with $n_i\in \mathbb{N}$ for $1 \leq i \leq k$. We write $(\widehat{S}, \widehat{g_i^G})=n_i$ following \cite{AradFisman}. Nevertheless, our notation for class sums differs from that appearing in \cite{AradFisman} in order to facilitate the reading.\\

\noindent
\textbf{Lemma 2.3.} {\it If $D_{1}$, $D_{2}$ and $D_{3}$ are conjugacy classes of a finite group $G$, then}
\begin{enumerate}[label=(\roman*)]
\item $(\widehat{D_{1}}\widehat{D_{2}}, \widehat{D_{3}})=(\widehat{D_{1}^{-1}}\widehat{D_{2}^{-1}}, \widehat{D_{3}^{-1}})$
\item $(\widehat{D_{1}}\widehat{D_{2}}, \widehat{D_{3}})=|D_{2}||D_{3}|^{-1}(\widehat{D_{1}}\widehat{D_{3}^{-1}}, \widehat{D_{2}^{-1}})$
\item $(\widehat{D_{1}}\widehat{D_{2}}, \widehat{D_{1}})=|D_{2}||D_{1}|^{-1}(\widehat{D_{1}}\widehat{D_{1}^{-1}}, \widehat{D_{2}^{-1}})=(\widehat{D_{2}}\widehat{D_{1}^{-1}}, \widehat{D_{1}^{-1}})=(\widehat{D_{2}^{-1}}\widehat{D_{1}}, \widehat{D_{1}})$.
\end{enumerate}

{\it Proof.} See the proof of Theorem A of \cite{AradFisman}. $\Box$

\section{Proofs}

We start by proving Theorem C, which is used for proving Theorem B.\\

{\it Proof of Theorem C.} The case $A=A^{-1}$ is easy and known, so we can assume that $A\neq A^{-1}$. By Lemma 2.3 we have
$$m=(\widehat{A}\widehat{A^{-1}}, \widehat{A})=(\widehat{A}\widehat{A^{-1}}, \widehat{A^{-1}})=(\widehat{A^2}, \widehat{A})=(\widehat{A^{-2}}, \widehat{A^{-1}}),$$

\noindent
where $m\in \mathbb{N^+}$ and so we can write 
\begin{equation}
\begin{split}
\widehat{A}\widehat{A^{-1}}=|A|\widehat{1}+m\widehat{A}+m\widehat{A^{-1}}\\
\widehat{A^2}=m\widehat{A}+\alpha \widehat{A^{-1}}+\widehat{T}\\
 \widehat{A^{-2}}=m\widehat{A^{-1}}+\alpha \widehat{A}+\widehat{T^{-1}}
 \end{split}
 \end{equation}
 
\noindent
with $\alpha\in \mathbb{N^+}$ and $\widehat{T}$ is a sum of conjugacy classes such that $(\widehat{T}, \widehat{L})=(\widehat{T^{-1}}, \widehat{L})=0$ for $L\in \{1, A, A^{-1}\}$.\\
 
 Suppose $T\neq \emptyset$ and calculate
 $$\widehat{A^2}\widehat{A^{-2}}=(m\widehat{A}+\alpha \widehat{A^{-1}}+\widehat{T})(m\widehat{A^{-1}}+\alpha \widehat{A}+\widehat{T^{-1}})=$$
 $$=m^2\widehat{A}\widehat{A^{-1}}+m\alpha \widehat{A^2}+m\widehat{A}\widehat{T^{-1}}+m\alpha \widehat{A^{-2}}+\alpha^2 \widehat{A^{-1}}\widehat{A}+\alpha \widehat{A^{-1}}\widehat{T^{-1}}+m\widehat{T}\widehat{A^{-1}}+\alpha \widehat{T}\widehat{A}+\widehat{T}\widehat{T^{-1}}.$$
 
\noindent 
For convenience, we write $$\widehat{T}=l_1\widehat{L_1}+\cdots+l_s\widehat{L_s}$$
where $L_i$ are distinct conjugacy classes of $G$ and $l_i$ the corresponding multiplicities. Consequently, from the above equation, we observe
 
 \begin{equation}
 (\widehat{A^2}\widehat{A^{-2}}, \widehat{1})=m^2|A|+\alpha^2|A|+l_1|L_1|+\cdots l_s|L_s|.
 \end{equation}
 
\noindent 
On the other hand, 

$$(\widehat{A}\widehat{A^{-1}})^2=(|A|\widehat{1}+m\widehat{A}+m\widehat{A^{-1}})(|A|\widehat{1}+m\widehat{A}+m\widehat{A^{-1}})=$$
$$=|A|^2\widehat{1}+|A|m\widehat{A}+|A|m\widehat{A^{-1}}+|A|m\widehat{A}+m^2\widehat{A^2}+m^2\widehat{A}\widehat{A^{-1}}+m|A|\widehat{A^{-1}}+m^2\widehat{A^{-1}}\widehat{A}+m^2\widehat{A^{-2}}.$$

\noindent
Thus, 
 \begin{equation}
 ((\widehat{A}\widehat{A^{-1}})^2, \widehat{1})=|A|^2+2m^2|A|.
 \end{equation}
\noindent
By joining Eqs. (2) and (3) we obtain
\begin{equation}
l_1|L_1|+\cdots l_s|L_s|=|A|^2+(m^2-\alpha^2)|A|
\end{equation}
\noindent
and from Eq. (1) we have 
\begin{equation}
l_1|L_1|+\cdots l_s|L_s|=|A|^2-(m+\alpha)|A|.
\end{equation}
\noindent
Hence, from Eqs. (4) and (5), we conclude that $m=\alpha-1$.\\

On the other hand, we calculate 
\begin{equation}
\begin{split}
\widehat{A}(\widehat{A}\widehat{A^{-1}})=\widehat{A}(|A|\widehat{1}+m\widehat{A}+m\widehat{A^{-1}})=\\
=|A|\widehat{A}+m(m\widehat{A}+\alpha \widehat{A^{-1}}+\widehat{T})+m(|A|\widehat{1}+m\widehat{A}+m\widehat{A^{-1}})=\\
=m|A|\widehat{1}+(|A|+2m^2)\widehat{A}+(\alpha m+m^2)\widehat{A^{-1}}+m\widehat{T}
\end{split}
\end{equation}
\noindent
and 
\begin{equation}
\begin{split}
(\widehat{A^2})\widehat{A^{-1}}=(m\widehat{A}+\alpha \widehat{A^{-1}}+\widehat{T})\widehat{A^{-1}}=\\=m(|A|\widehat{1}+m\widehat{A}+m\widehat{A^{-1}})+\alpha(m\widehat{A^{-1}}+\alpha \widehat{A}+\widehat{T^{-1}})+\widehat{T}\widehat{A^{-1}}=\\
=|A|m\widehat{1}+(m^2+\alpha^2)\widehat{A}+(m^2+\alpha m)\widehat{A^{-1}}+\alpha \widehat{T^{-1}}+\widehat{T}\widehat{A^{-1}}. 
\end{split}
\end{equation}

So, from Eqs. (6) and (7) we conclude that $T=T^{-1}$ and $m=\alpha+\beta$ for some $\beta \in \mathbb{N}^*$, a contradiction. As a consequence, $T=\emptyset$ and $A^2=A \cup A^{-1}$.\\

Now we prove that $\langle A\rangle$ is elementary abelian. Indeed, we have $A^3=AA^2=A(A\cup A^{-1})=A^2 \cup AA^{-1}=1 \cup A \cup A^{-1}$, so we deduce that $\langle A\rangle=1\cup A \cup A^{-1}$. In particular, all non-trivial elements of $\langle A\rangle$ have the same order, and this forces $\langle A\rangle$ to be $p$-elementary for some prime $p$. Finally, we prove that $\langle A\rangle$ is abelian. Put $N=\langle A\rangle$ and let $x\in A$. Observe that $|x^N|$ divides $|A|=|x^G|$, but on the other hand, $|x^N|$ also divides $|N|=1+2|A|$. This implies that $|x^N|=1$, and hence $N$ is abelian. $\Box$\\

\textit{Examples.} The smallest group for Theorem C with $A$ non-trivial and real is the symmetric group on 3 letters with the conjugacy  class of 3-cycles.  The smallest example for Theorem C with $A$ non-real is the non-abelian group of order 21,  $G=\langle x, y\, \mid \, x^y =x^2, x^7=1\rangle$, when we consider the conjugacy class $A=\{x, x^2, x^4\}$ where $\langle A\rangle=\langle x\rangle\cong\mathbb{Z}_7$.\\

We restate Theorems A and B in terms of the theorems proved in \cite{AradFisman}. We will divide the proofs into steps, some of them being equal to those appearing in the original proofs in \cite{AradFisman}.\\ 

\noindent
\textbf{Theorem A.} {\it Let $D_1$ and $D_2$ be conjugacy classes of a finite group $G$ and suppose that $D_1D_2=D_1 \cup D_2$. Then $\langle D_1\rangle=\langle D_2\rangle$ is solvable. Furthermore, the elements in $D_1$ and $D_2$ are $p$-elements for some prime $p$ and $\langle D_1\rangle$ is $p$-nilpotent. }\\

{\it Proof.} First, let us prove $\langle D_1\rangle=\langle D_2\rangle$ by induction on $|G|$. If $G=\langle D_1 \rangle=\langle D_2 \rangle$ the proof is finished. Suppose, for instance, that $\langle D_1 \rangle < G$ and write $\overline{G}=G/\langle D_1 \rangle$. Then $\overline{D_1}\overline{D_2}=\overline{D_1}\cup \overline{D_2}$, which implies that $\overline{D_2}=\overline{1}$, and hence $\langle D_2 \rangle \subseteq \langle D_1 \rangle$. Now we consider $\widehat{G}=G/\langle D_2 \rangle$ and, arguing as above, $\langle D_1 \rangle \subseteq \langle D_2 \rangle$.\\

We continue the proof by induction on $|G|$. We write $\widehat{D_1}\widehat{D_2}=n_1\widehat{D_1}+n_2\widehat{D_2}$ with $n_1, n_2 \in \mathbb{N}^*$.\\

\textbf{Step 1:} \textit{$\widehat{D_1}\widehat{D_2^{-1}}=n_1\widehat{D_1}+n_2\widehat{D_2^{-1}}$ and $D_i=D_{i}^{-1}$ for $1 \leq i \leq 2$.} \\

Mimic steps b(i) and b(ii) of the proof of Theorem 2 of \cite{AradFisman}.\\

\textbf{Step 2:} \textit{We have} $$\widehat{D_{1}^{2}}=|D_1|\widehat{1}+n_1|D_1||D_2|^{-1}\widehat{D_2}+s_1\widehat{D_1}+\widehat{M_1}$$ $$\widehat{D_{2}^{2}}=|D_2|\widehat{1}+n_2|D_2||D_1|^{-1}\widehat{D_1}+s_2\widehat{D_2}+\widehat{M_2}$$
\textit{where} $s_i\in \mathbb{N}$ \textit{and} $\widehat{M_i}$ are sums of conjugacy classes such that $(\widehat{M_i}, \widehat{C})=0$ \textit{for} $C \in \{1, D_j\}$, $i, j, \in \{1, 2\}.$\\

By Lemma 2.3 we know  that $$(\widehat{D_1^2}, \widehat{1})=|D_1|(\widehat{D_1}\widehat{1}, \widehat{D_1})=|D_1|,$$
$$(\widehat{D_1^{2}}, \widehat{D_2})=|D_1||D_2|^{-1}(\widehat{D_1}\widehat{D_2}, \widehat{D_1})=|D_1||D_2|^{-1}n_1.$$

\noindent
Then we can write $$\widehat{D_{1}^{2}}=|D_1|\widehat{1}+n_1|D_1||D_2|^{-1}\widehat{D_2}+s_1\widehat{D_1}+\widehat{M_1}$$ and analogously, $$\widehat{D_{2}^{2}}=|D_2|\widehat{1}+n_2|D_2||D_1|^{-1}\widehat{D_1}+s_2\widehat{D_2}+\widehat{M_2}$$
for some $s_i\in \mathbb{N}$ and $\widehat{M_i}$ such that $(\widehat{M_i}, \widehat{C})=0$ for $C \in \{1, D_j\}$, $i, j, \in \{1, 2\}.$\\

We distinguish two subcases depending on whether $M_1=\emptyset$ or not.\\

\textbf{Step 3:} \textit{If} $M_1=\emptyset$\textit{, then} $\langle D_1\rangle$ \textit{is} $p$\textit{-elementary abelian, so the theorem is proved.}\\

If $M_1=\emptyset$, then either $D_1^2=1 \cup D_1 \cup D_2$ or $D_1^2= 1 \cup D_2$. In the first case, $D_1^3=1\cup D_1 \cup D_2$ and in the second $D_1^4=1\cup D_1 \cup D_2$, so in both cases it certainly follows that $\langle D_1\rangle=1 \cup D_1 \cup D_2$. Hence, joint with the fact that $\langle D_1\rangle=\langle D_2\rangle$, we deduce that $\langle D_1\rangle$ is a minimal normal subgroup of $G$. Furthermore, it must be solvable due to the fact that its elements only have two possible orders. Consequently, $\langle D_1\rangle$ is $p$-elementary abelian for some prime $p$, so the thesis of the theorem trivially follows. \\

Henceforth, we will assume that $M_1\neq \emptyset$.\\

\textbf{Step 4:} \textit{We have} 
$$n_1\widehat{M_1}=n_1|D_1||D_2|^{-1}\widehat{M_2}+\widehat{M_1}\widehat{D_2}-(\widehat{M_1}\widehat{D_2}, \widehat{D_2})\widehat{D_2}$$
$$n_2\widehat{M_2}=n_2|D_2||D_1|^{-1}\widehat{M_1}+\widehat{M_2}\widehat{D_1}-(\widehat{M_2}\widehat{D_1}, \widehat{D_1})\widehat{D_1}.$$

Mimic step b(iv) of the proof of Theorem 2 of \cite{AradFisman}.\\

\textbf{Step 5:} \textit{Conclusion.}\\

First, let us see that $\langle D_1\rangle$ is solvable. By applying step 4, we have
$$n_1n_2\widehat{M_2}=n_1n_2|D_2||D_1|^{-1}\widehat{M_1}+n_1(\widehat{M_2}\widehat{D_1}-(\widehat{M_2}\widehat{D_1}, \widehat{D_1})\widehat{D_1})=$$
$$=n_2|D_2||D_1|^{-1}(n_1|D_1||D_2|^{-1}\widehat{M_2}+\widehat{M_1}\widehat{D_2}-(\widehat{M_1}\widehat{D_2}, \widehat{D_2})\widehat{D_2})+$$
$$n_1\widehat{M_2}\widehat{D_1}-n_1(\widehat{M_2}\widehat{D_1}, \widehat{D_1})\widehat{D_1}.$$

\noindent
It follows that $n_2|D_2||D_1|^{-1}\widehat{M_1}\widehat{D_2}+n_1\widehat{M_2}\widehat{D_1}=l_1\widehat{D_1}+l_2\widehat{D_2}$ for $\{l_1, l_2\} \subset \mathbb{N}^*$. In particular, $\widehat{M_1}\widehat{D_2}=m_1\widehat{D_1}+m_2\widehat{D_2}$ for $\{m_1, m_2\} \subset \mathbb{N}^*$. Since $$(\widehat{M_1}\widehat{D_2}, \widehat{D_1})=\dfrac{|M_1|}{|D_1|}(\widehat{M_1}, \widehat{D_1}\widehat{D_2})=0,$$ then $\widehat{M_1}\widehat{D_2}=m_2\widehat{D_2}$. Symmetrically, $\widehat{M_2}\widehat{D_1}=m_1\widehat{D_1}$. Hence, taking into account step 4 and the definition of $M_1$, we have $n_1\widehat{M_1}=n_1|D_1||D_2|^{-1}\widehat{M_2}$. Thus, $M_1=M_2$, as a set that is union of conjugacy classes, so we have $D_1M_1=D_1$. As a result,  $D_1\langle M_1\rangle=D_1$. Then, as $D_1^2=1 \cup D_2 \cup D_1 \cup M_1$, we conclude that $\langle D_1\rangle=D_1^2=1 \cup D_2 \cup D_1 \cup M_1$.\\

We write $\overline{G}=G/\langle M_1\rangle$ and then $\langle \overline{D_1}\rangle=\overline{1} \cup \overline{D_1}\cup \overline{D_2}$. By induction, the elements in $\overline{D_1}$ and $\overline{D_2}$ are $p$-elements for some prime $p$, so $\langle \overline{D_1}\rangle$ is a $p$-group. Let $d\in D_1$. As all elements in $d\langle M_1\rangle$ are conjugate in $G$, then $\langle M_1\rangle$ is solvable by Theorem 2.1(a). It clearly follows that $\langle D_1\rangle$ is solvable. Finally let us prove that the elements in $D_1$ and $D_2$ are $p$-elements for some prime $p$. Let $P\in {\rm Syl}_{p}(\langle D_1\rangle)$. Note that  $\langle D_1 \rangle=P\langle M_1\rangle=P\langle M_1\rangle D_1=PD_1$. In particular, we can write $1=xd$ with $x\in P$ and $d\in D_1$. This shows that the elements in $D_1$ are $p$-elements. Analogously, we can deduce that the elements in $D_2$ are also $p$-elements. By Theorem 2.1(b), we conclude that $\langle M_1\rangle$, and then also $\langle D_1\rangle$, has normal $p$-complement. $\Box$\\

\textit{Examples.} We show different examples corresponding to distinct cases of  Theorem A. \\

1. The easiest example is the dihedral group of order 10, in which the only two conjugacy classes of size 2 satisfy the hypotheses of the theorem. This example can be generalized by taking  $G=\langle x\rangle \rtimes \langle a\rangle \cong \mathbb{Z}_p \rtimes \mathbb{Z}_{(p-1)/2}$, where $p$ is a prime such that $p \equiv 1\, ({\rm mod}\, 4)$ and $a$ is an automorphism of order $(p-1)/2$ of $\langle x\rangle$.  The subgroup $\langle x\rangle$ contains exactly the trivial class and two (real) conjugacy classes $A$ and $B$ of size $(p-1)/2$, which satisfy $AB= A \cup B$ and  $\langle A\rangle=\langle x\rangle$. This corresponds to the case $M_1=\emptyset$ in step 3 of the proof of Theorem A.\\

2. Two examples with $\langle A\rangle$ non-cyclic are the following. Let $G=(\langle x\rangle \times \langle y\rangle) \rtimes (\langle a\rangle \times \langle b\rangle)  \cong (\mathbb{Z}_3 \times \mathbb{Z}_3) \rtimes (\mathbb{Z}_2\times \mathbb{Z}_2)$, where $a$ and $b$ are defined by: $x^a=x^2$, $y^a=y^2$, $x^b=y$ and  $y^b=x$. When we take $A=\{x,x^2, y,y^2\}$ and $B=\{xy, x^2y^2, xy^2, x^2y\}$ then we have $AB= A \cup B$ and $\langle A\rangle=\langle B\rangle=\langle x\rangle \times \langle y\rangle$. On the other hand,  the group of the library of the small groups of GAP \cite{GAP} with number ${\rm Id}(1176, 213)$ has two conjugacy classes $A$ and $B$ of size 24 satisfying the hypotheses of Theorem A, with $\langle A\rangle\cong \mathbb{Z}_7\times \mathbb{Z}_7$. Also in this case $M_1=\emptyset$.\\

3. The group ${\rm Id}(108, 15)$ has two conjugacy classes $A$ and $B$ of size 12 satisfying $AB=A\cup B$ with $\langle A\rangle\cong (\mathbb{Z}_3\times \mathbb{Z}_3) \rtimes \mathbb{Z}_3$. This example shows that $\langle A\rangle$ is not necessarily abelian. We remark that this example corresponds to the case $M_1\neq \emptyset$ in the proof of Theorem A (see step 4). In fact, $\langle M_1\rangle=\textbf{Z}(\langle A\rangle)$.\\
 
 As we said in the Introduction, to prove Theorem B we make a slight improvement of Theorem D of \cite{Powers} by adding $p$-nilpotency.\\
 
 \noindent
\textbf{Theorem 3.1.} {\it  Let $G$ be a group and let $K = x^G$ be a conjugacy class of $G$. If $K^2 = K \cup K^{-1}$, then $\langle K\rangle$ is solvable. Moreover, $x$ is a p-element for some prime $p$ and $\langle K\rangle$ is $p$-nilpotent.}\\

{\it Proof.} Following the proof of Theorem D of \cite{Powers} we have $KK^{-1}=1 \cup K \cup K^{-1} \cup S$ where $S$ is union of conjugacy classes of $G$ other than $1$, $K$  and $K^{-1}$. If $S=\emptyset$, by the original proof, then $\langle K\rangle$ is $p$-elementary abelian for some prime $p$ and the theorem is proved. If $S\neq \emptyset$, again by following the original proof, $KS=K$, $\langle K \rangle/\langle S \rangle$ is $p$-elementary abelian for some prime $p$ and $x$ is a $p$-element. In particular, the elements of $x\langle S\rangle$ are all conjugate in $G$ and, by applying Theorem 2.1(b), $\langle S\rangle$ has normal $p$-complement. Since $\langle K \rangle/\langle S \rangle$ is a $p$-group, then $\langle K\rangle$ has normal $p$-complement too. $\Box$
\bigskip

\noindent
\textbf{Theorem B.} {\it Let $D_1$ and $D_2$ be conjugacy classes of a finite group $G$ and suppose that $D_1D_2=D_1^{-1} \cup D_2$ with $D_1\neq D_1^{-1}$. Then $D_1=D_2$ and $\langle D_1\rangle$ is solvable. Moreover, $D_1$ is a class of $p$-elements and $\langle D_1\rangle$ is  $p$-nilpotent.}\\

{\it Proof.} By arguing by induction on $|G|$ as at the beginning of the proof of Theorem A, we can easily deduce that $\langle D_1\rangle=\langle D_2\rangle$.\\

If $D_1 =D_2$, by Theorem 3.1 we have that $\langle D_1\rangle$ is solvable, the elements of $D_1$ are $p$-elements, and $\langle D_1\rangle$ is $p$-nilpotent. To complete the proof, in the following, we will prove by minimal counterexample that there do not exist distinct classes $D_1$ and $D_2$ in a finite group satisfying the hypotheses of the theorem. Let $G$ be a finite group of minimal order and let $D_1$ and $D_2$ two conjugacy classes such that $D_1D_2= D_1^{-1} \cup D_2$, with $D_1$ non-real and $D_1\neq D_2$. We write $\widehat{D_1}\widehat{D_2}=n_1\widehat{D_1^{-1}}+n_2\widehat{D_2}$ with $n_1, n_2 \in \mathbb{N}^*$. We distinguish two cases: first, $D_2=D_2^{-1}$ and second $D_2\neq D_2^{-1}$.\\

\textbf{Case 1:} $D_2=D_2^{-1}$.\\

\textbf{Step 1.1:} \textit{We have $$\widehat{D_{1}^{2}}=n_1\widehat{D_2}+n_2\widehat{D_1}$$
$$\widehat{D_{2}^{2}}=|D_2|\widehat{1}+n_2(\widehat{D_1}+\widehat{D_1^{-1}})+\widehat{L}$$ with $L=L^{-1}$, $0=(\widehat{L}, \widehat{C})$ for $C \in \{1, D_1, D_1^{-1}, D_2\}$.} \\

Follow steps c(1)(i), c(1)(ii), c(1)(iii) and c(1)(iv) of the proof of Theorem 2 of \cite{AradFisman}.\\

\textbf{Step 1.2:} \textit{We have} $\widehat{D_1^{-1}}\widehat{D_2}=n_1\widehat{D_1}+n_2\widehat{D_2}$.\\

Since $\widehat{D_1}\widehat{D_2}=n_1\widehat{D_1^{-1}}+n_2\widehat{D_2}$ and $\widehat{D_2}=\widehat{D_2^{-1}}$, we have $\widehat{D_1^{-1}}\widehat{D_2}=\widehat{D_1^{-1}}\widehat{D_2^{-1}}=n_1\widehat{D_{1}}+n_2\widehat{D_{2}^{-1}}=n_1\widehat{D_{1}}+n_2\widehat{D_{2}}$.\\

\textbf{Step 1.3:} $L\neq \emptyset$.\\ 

If $L=\emptyset$, then $D_2^2=1 \cup D_1 \cup D_1^{-1}$ and since $D_2$ is real, by Lemma 2.2, $D_1$ is also real, again a contradiction.\\

\textbf{Step 1.4:} \textit{Conclusion.} \\

We know, by step 1.1, 
$$\widehat{D_2^2}\widehat{D_1}=(|D_2|\widehat{1}+n_2(\widehat{D_1}+\widehat{D_1^{-1}})+\widehat{L})\widehat{D_1}$$
$$=|D_2|\widehat{D_1}+n_2\widehat{D_1^2}+n_2\widehat{D_2^2}+\widehat{L}\widehat{D_1}$$ and, by applying step 1.2, $$\widehat{D_2}(\widehat{D_1}\widehat{D_2})=\widehat{D_2}(n_1\widehat{D_1^{-1}}+n_2\widehat{D_2})=n_1(n_1\widehat{D_1}+n_2\widehat{D_2})+n_2\widehat{D_2^2}.$$
Hence $$|D_2|\widehat{D_1}+n_2(n_1\widehat{D_2}+n_2\widehat{D_1})+\widehat{L}\widehat{D_1}=n_1^2\widehat{D_1}+n_1n_2\widehat{D_2}.$$

\noindent
Thus $(|D_2|+n_2^2)\widehat{D_1}+\widehat{L}\widehat{D_1}=n_1^2\widehat{D_1}$. It follows that $\widehat{L}\widehat{D_1}=k\widehat{D_1}$ for some $k \in \mathbb{N^*}$. As a consequence, there exists a conjugacy class $C$ of $G$ other than $1$, $D_1$, $D_1^{-1}$ and $D_2$ such that $D_1C=D_1$.  Thus, $D_1\langle C\rangle=D_1$, with $1\neq \langle C\rangle\triangleleft G$ and we write $\overline{G}=G/\langle C\rangle$. We have $|\overline{G}|< |G|$, $\overline{D_1}\overline{D_2}=\overline{D_1^{-1}}\cup \overline{D_2}$, with $\overline{D_1}\neq \overline{D_1^{-1}}$, because otherwise $D_1=D_1\langle C\rangle=D_1^{-1}\langle C\rangle=D_1^{-1}$, a contradiction. In addition, $\overline{D_1}\neq \overline{D_2}$ because otherwise $D_1=D_1\langle C\rangle=D_2\langle C\rangle\supseteq D_2$, which is impossible. By minimality, we get a contradiction and this case is finished. \\

\textbf{Case 2:} $D_2\neq D_2^{-1}$.\\

We have $$0=(\widehat{D_1}\widehat{D_2}, \widehat{D_1})=\dfrac{|D_2|}{|D_1|}(\widehat{D_1}\widehat{D_{1}^{-1}}, \widehat{D_{2}^{-1}})=\dfrac{|D_2|}{|D_1|}(\widehat{D_1}\widehat{D_{1}^{-1}}, \widehat{D_{2}})=(\widehat{D_1}\widehat{D_{2}^{-1}}, \widehat{D_1})$$

$$0=(\widehat{D_1}\widehat{D_2}, \widehat{D_{2}^{-1}})=\dfrac{|D_1|}{|D_2|}(\widehat{D_{2}^2}, \widehat{D_1^{-1}})$$

$$n_1=(\widehat{D_1}\widehat{D_2}, \widehat{D_{1}^{-1}})=\dfrac{|D_2|}{|D_1|}(\widehat{D_{1}^2}, \widehat{D_2^{-1}})$$

$$n_2=(\widehat{D_1}\widehat{D_2}, \widehat{D_2})=(\widehat{D_1}\widehat{D_{2}^{-1}}, D_2^{-1})=\dfrac{|D_1|}{|D_2|}(\widehat{D_2}\widehat{D_{2}^{-1}}, \widehat{D_{1}^{-1}}).$$
We denote by 
$$l_1=(\widehat{D_1}\widehat{D_2^{-1}}, \widehat{D_{1}^{-1}})=\dfrac{|D_2|}{|D_1|}(\widehat{D_1^2}, \widehat{D_2}), \, \, \, \, l_2=(\widehat{D_1}\widehat{D_2^{-1}}, \widehat{D_2})=\dfrac{|D_1|}{|D_2|}(\widehat{D_2^2}, \widehat{D_1})$$
$$j_1=(\widehat{D_1^2}, \widehat{D_1})=(\widehat{D_1}\widehat{D_1^{-1}}, \widehat{D_1}), \, \, \, \, j_2=(\widehat{D_2^2}, \widehat{D_2})=(\widehat{D_2}\widehat{D_2^{-1}}, \widehat{D_2^{-1}})$$
$$d_1=(\widehat{D_1^2}, D_1^{-1}), \, \, \, \, d_2=(\widehat{D_2^2}, D_2^{-1})$$
Therefore, we can collect all these multiplicities in Table 1, which also appears in the original proof.\\

\begin{table}
\centering
\begin{tabular}{c|cccccc}
 & $\widehat{1}$ & $\widehat{D_1}$ & $\widehat{D_1^{-1}}$ & $\widehat{D_2}$ & $\widehat{D_2^{-1}}$ & \\ \hline
 $\widehat{D_1}\widehat{D_2}$ & 0 & 0 & $n_1$ & $n_2$ & 0 & \\
  $\widehat{D_1}\widehat{D_2^{-1}}$ & 0 & 0 & $l_1$ & $l_2$ & $n_2$ & $\widehat{M_{12}}$\\
  $\widehat{D_1^2}$ & 0 & $j_1$ & $d_1$ & $l_1\dfrac{|D_1|}{|D_2|}$ & $n_1\dfrac{|D_1|}{|D_2|}$ & $\widehat{M_{11}}$\\
    $\widehat{D_2^2}$ & 0 & $l_2\dfrac{|D_2|}{|D_1|}$  & 0 & $j_2$ & $d_2$ & $\widehat{M_{22}}$\\
    $\widehat{D_1}\widehat{D_{1}^{-1}}$ & $|D_1|$ & $j_1$ & $j_1$ & 0 & 0 & $N_1$\\
    $\widehat{D_2}\widehat{D_{2}^{-1}}$ & $|D_2|$ & $n_2\dfrac{|D_2|}{|D_1|}$ & $n_2\dfrac{|D_2|}{|D_1|}$ & $j_2$ & $j_2$ & $N_2$
\end{tabular}
\caption{}
\label{T1}
\end{table}

\noindent
In Table 1, we have $N_i=N_{i}^{-1}$ and $(\widehat{L}, \widehat{C})=0$ for $C\in \{1, D_k, D_{k}^{-1}\}$, $L\in \{M_{ij}, N_i\}$ for every $k, i, j \in \{1, 2\}$.\\

\textbf{Step 2.1:} $n_1\widehat{N_1}=n_1\dfrac{|D_1|}{|D_2|}\widehat{N_2}$ \textit{and} $\widehat{N_2}\widehat{D_1}=(\widehat{N_2}\widehat{D_1}, \widehat{D_1})\widehat{D_1}$.\\

Follow steps c(2)(i) to (vii) of the proof of Theorem 2 of \cite{AradFisman}.\\

\textbf{Step 2.2:} \textit{Conclusion.}\\

We distinguish two cases, whether $\widehat{N_2}\neq 0$ or not. First, if $\widehat{N_2}\neq 0$, then there exists a conjugacy class $C$ of $G$ such that $D_1C=D_1$. We can apply the same argument as at the end of step 1.4 of Case 1 and, by minimal counterexample, we get a contradiction.\\

Assume now that $\widehat{N_2}=0$. By step 2.1, we know that $$n_1\widehat{N_1}=n_1\dfrac{|D_1|}{|D_2|}\widehat{N_2}.$$ Therefore, $\widehat{N_1}=0$. Thus, from Table 1, we have $D_1D_1^{-1}= 1 \cup D_1 \cup D_1^{-1}$ and by Theorem C, we conclude that $\langle D_1\rangle=D_1D_1^{-1}= 1 \cup D_1 \cup D_1^{-1}$. This forces that $D_2=D_1^{-1}$ or $D_2=D_1$ and both certainly are contradictions. This finishes the proof. $\Box$\\

\textit{Example.} This is an example of Theorem B where $\langle K \rangle$ is $p$-nilpotent and not a $p$-group. We take the group $G=((\mathbb{Z}_2 \times \mathbb{Z}_2 \times \mathbb{Z}_2) \rtimes \mathbb{Z}_7) \rtimes \mathbb{Z}_3={\rm Id}(168,43)$ which has a conjugacy class $K$ of elements of order 7 and size 24 satisfying $K^2=K \cup K^{-1}$. Also, $\langle K\rangle=(\mathbb{Z}_2 \times \mathbb{Z}_2 \times \mathbb{Z}_2) \rtimes \mathbb{Z}_7$.\\
 
\noindent {\bf Acknowledgements.} This research is partially supported by the Spanish Government, Proyecto PGC2018-096872-B-I00. The first author is also supported by Proyecto UJI-B2019-03.

\end{document}